\documentclass{amsart}

\usepackage{amssymb}
\usepackage[english]{babel}

\newtheorem{theorem}{Theorem}[section]
\newtheorem{lemma}[theorem]{Lemma}
\newtheorem{e-proposition}[theorem]{Proposition}
\newtheorem{corollary}[theorem]{Corollary}
\newtheorem{e-definition}[theorem]{Definition\rm}
\newtheorem{remark}{\it Remark\/}

\newtheorem{theoreme}{Th\'eor\`eme}[section]

\newtheorem{proposition}[theoreme]{Proposition}

\newcommand{\RR}{\mathbb{R}}

\newcommand{\cD}{{{\mathcal{D}}}}
\newcommand{\cM}{{{\mathcal{M}}}}
\newcommand{\supp}{\operatorname{supp}}
\newcommand{\capt}{\operatorname{cap}}
\newcommand{\sign}{\operatorname{sgn}}
\newcommand{\loc}{_{\rm loc}}

\newcommand{\nd}{\noindent}
\newcommand{\eps}{\varepsilon}
\newcommand{\lc}{_{\rm c}}
\newcommand{\ld}{_{\rm d}}

\def\og{\leavevmode\raise.3ex\hbox{$\scriptscriptstyle\langle\!\langle$~}}
\def\fg{\leavevmode\raise.3ex\hbox{~$\!\scriptscriptstyle\,\rangle\!\rangle$}}

\begin{document}

\title{Kato's inequality when $\Delta u$ is a measure}

\author{Ha\"\i m Brezis}
\address{
Ha\"\i m Brezis\hfill\break\indent
Universit{\'e} Pierre et Marie Curie\hfill\break\indent
Laboratoire Jacques-Louis Lions\hfill\break\indent
4 pl.\@ Jussieu, BC\@ 187\hfill\break\indent
75252 Paris Cedex 05, France\hfill\break\indent
\medskip
Rutgers University\hfill\break\indent
Dept.\@ of Math., Hill Center, Busch Campus\hfill\break\indent
110 Frelinghuysen Rd\hfill\break\indent
Piscataway, NJ 08854, USA}

\author{Augusto C. Ponce}
\address{
Augusto C. Ponce\hfill\break\indent
Institute for Advanced Study\hfill\break\indent
Princeton, NJ 08540, USA
}

\begin{abstract}
% Text of abstract in English
We extend the classical Kato's inequality in order to allow functions
 $u \in L^1\loc$ such that $\Delta u$ is a Radon measure. This inequality
 has been recently applied by Brezis, Marcus, and
 Ponce~\cite{BreMarPon:XX} to study the existence of solutions
 of the nonlinear equation $- \Delta u + g(u) = \mu$, where $\mu$ is a measure
 and $g : \RR \to \RR$ is an increasing continuous function.\\
 To cite this paper: Ha\"\i m~Brezis and Augusto C.~Ponce, \textit{Kato's inequality 
  when $\Delta u$ is a measure}, C.\@ R.\@ Acad.\@ Sci.\@ Paris,
  Ser.\@ I  \textbf{338} (2004), 599--604.
\end{abstract}

\maketitle

\section{Introduction and main result}
\label{sec1}

Let $N \geq 1$ and $\Omega \subset \RR^N$ be a bounded open subset.
The classical Kato's inequality (see \cite{Kato:72}) states that given
any function $u \in L^1\loc(\Omega)$ such that $\Delta u \in L^1\loc(\Omega)$,
then $\Delta u^+$ is a Radon measure and the following holds:
\begin{equation}\label{1.1}
\Delta u^+ \geq \chi_{[u \geq 0]} \Delta u \quad \mbox{in $\cD^\prime(\Omega)$.}
\end{equation}

Our main result in this paper (see Theorem~\ref{thm1.1} below) extends
(\ref{1.1}) to the case 
$\Delta u \in \cM(\Omega)$, where $\cM(\Omega)$ denotes the space of 
Radon measures on $\Omega$. In other words, $\mu \in \cM(\Omega)$ if
and only if, for every $\omega \subset\subset \Omega$, there exists
$C_\omega>0$ such that 
$\left| \int_\Omega \varphi \, d\mu \right| \leq C_\omega \|\varphi\|_\infty$,
$\forall \varphi \in C_0^\infty(\omega)$.

We first recall that any $\mu \in \cM(\Omega)$ can be uniquely decomposed as a sum
of two Radon measures on $\Omega$ (see e.g. \cite{FukSatTan:91}):
$\mu = \mu\ld + \mu\lc$,
where
\begin{eqnarray*}
\mu\ld(A) &=& 0 \quad \mbox{for any Borel measurable set $A \subset
  \Omega$ such that $\capt{(A)} = 0$,}\\
|\mu\lc|(\Omega \backslash F)  &=& 0 \quad \mbox{for some Borel
  measurable set $F \subset \Omega$ such that $\capt{(F)} = 0$.}
\end{eqnarray*}
Here, $\capt$ denotes the Newtonian ($W^{1,2}$) capacity of a set. We observe
that $\mu\ld$ and $\mu\lc$ are singular with respect to each other.
This decomposition is the analog of the classical Radon-Nikodym
Theorem, but with respect to $\capt$. Clearly, $(\mu\ld)^+ =
(\mu^+)\ld$ and $(\mu\lc)^+ = (\mu^+)\lc$.

Using the above notation, we can now state our main result:

\smallskip
\begin{theorem}\label{thm1.1}
Let $u \in L^1\loc(\Omega)$ be such that $\Delta u \in \cM(\Omega)$. Then,
$\Delta u^+ \in \cM(\Omega)$, and the following holds:
\begin{alignat}{2}
(\Delta u^+)\ld  & \geq \chi_{[u \geq 0]} (\Delta u) \ld && \quad  \mbox{on $\Omega$,}\label{1.3}\\
(- \Delta u^+)\lc & = (- \Delta u)^+\lc \hspace{3.8ex} && \quad  \mbox{on $\Omega$.}\label{1.4}
\end{alignat}
\end{theorem}
\smallskip

Note that the right-hand side of (\ref{1.3}) is well-defined because
$u$ is quasicontinuous. More precisely, if $u \in L^1\loc(\Omega)$ and $\Delta u \in
\cM(\Omega)$, then there exists $\tilde u : \Omega \to \RR$ quasicontinuous
such that $u = \tilde u$ a.e. in $\Omega$ (see \cite{Anc:79} and also
\cite[Lemma~1]{BrePon:03}). In (\ref{1.3}), we then identify $u$ with
its quasicontinuous representative. It is easy to see that $\chi_{[u
    \geq 0]}$ is locally integrable in $\Omega$ with respect to the measure $\big|
(\Delta u)\ld \big|$.

The proof of (\ref{1.3}) requires a theorem of Boccardo, Gallou\"et,
and Orsina~\cite{BocGalOrs:96}, which says that a Radon measure
$\mu$ is diffuse (i.e. $\mu\lc = 0$) if and only if $\mu \in L^1\loc(\Omega) 
+ \Delta\big[H^1\loc(\Omega) \big]$. Identity (\ref{1.4}) relies on
(and in fact is equivalent to) the ``inverse'' maximum 
principle, recently established by Dupaigne and Ponce~\cite{DupPon:04}
(see Theorem~\ref{lemma3.1} below).

An equivalent statement of Theorem~\ref{thm1.1} is the following:

\smallskip
\begin{corollary}\label{thm1.2}
Let $u \in L^1\loc(\Omega)$ be such that $\Delta u \in \cM(\Omega)$. Then,
$\Delta |u| \in \cM(\Omega)$, and the following holds:
\begin{eqnarray}
(\Delta |u|)\ld  & \geq & \sign{(u)} \, (\Delta u) \ld \quad  \mbox{on
    $\Omega$,}\label{1.5}\\ 
(\Delta |u|)\lc & = & - |\Delta u|\lc \hspace{6.2ex} \quad  \mbox{on
    $\Omega$.}\label{1.6} 
\end{eqnarray}
\end{corollary}

\nd
Here, $\sign{(t)} = 1$ for $t>0$,  $\sign{(t)} = -1$ for $t<0$, and
$\sign{(0)} = 0$. 

\smallskip

\begin{remark}\label{rem1.1}
A slight modification of the proof of Theorem~\ref{thm1.1} shows that
\begin{equation}\label{1.7}
(\Delta u^+)\ld \geq \chi_{[u > 0]} (\Delta u)\ld \quad \mbox{on $\Omega$.}
\end{equation}
In other words, we can replace the set $[u \geq 0]$ in (\ref{1.3}) by
$[u > 0]$ and still get the same result.
\end{remark}
\smallskip

Here is a simple consequence of (\ref{1.7}):

\smallskip
\begin{corollary}\label{cor1.1}
Let $u \in L^1\loc(\Omega)$ be such that $\Delta u \in \cM(\Omega)$. If $u
\geq 0$ a.e. in $\Omega$, then
\begin{equation}\label{1.8}
(\Delta u)\ld \geq 0 \quad \mbox{on the set\/ $[u = 0]$.} 
\end{equation}
\end{corollary}

%%%%%%%%%%%%%%%%%%%%%%%%%%%%%%%%%%%%%%%%
%%%%%%%%%%%%%%%%%%%%%%%%%%%%%%%%%%%%%%%%
%%%%%%%%%%%%%%%%%%%%%%%%%%%%%%%%%%%%%%%%

\section{Proof of (\ref{1.3}) in Theorem~\ref{thm1.1}}

We start with the following:

\smallskip
\begin{lemma}\label{lemma2.1}
Assume $\mu \in \cM(\Omega)$ is a diffuse measure with respect to
$\capt$ (i.e.\@ $\mu\lc = 0$ on $\Omega$). Let $(v_n)$ be a
sequence in $L^\infty(\Omega) \cap H^1(\Omega)$ such that
$\|v_n\|_\infty \leq C$ and $v_n \rightharpoonup v$ in $H^1$.
Then,
\begin{equation}\label{2.1}
v_n \to v \quad \mbox{in $L^1\loc(\Omega; d\mu)$.}
\end{equation}
Equivalently, there exists a subsequence $(v_{n_k})$ converging to $v$ $|\mu|$-a.e.\@ in $\Omega$.
\end{lemma}
\smallskip

\noindent
\textbf{Proof.} Without loss of generality, we may assume that
$|\mu|(\Omega) < \infty$. By Theorem~2.1 of Boccardo, Gallou\"et, and
Orsina~\cite{BocGalOrs:96}, we know that 
$\mu = f - \Delta g$ in $\cD'(\Omega)$, for some $f \in
L^1(\Omega)$ and $g \in H^1(\Omega)$.
Using a standard density argument, we conclude that
\begin{equation}\label{2.2}
\int_\Omega w \varphi \, d\mu = \int_\Omega w \varphi f + \int_\Omega
\nabla g \cdot \nabla(w \varphi), \quad \forall \varphi \in C_0^\infty(\Omega), \quad
\forall w \in L^\infty \cap H^1.
\end{equation}
By assumption, the sequence $\big(|v_n - v| \big)$ is bounded in
$H^1(\Omega)$ and, by Rellich's theorem, $|v_n -v| \to 0$ in
$L^2(\Omega)$. Thus,
\begin{equation}\label{2.5}
|v_n - v| \rightharpoonup 0 \quad \mbox{in $H^1$.}
\end{equation}
Given $\eps > 0$, let $\omega \subset\subset \Omega$ be such that
$|\mu|(\Omega \backslash \omega) < \eps$. We then fix $\varphi_0 \in
C_0^\infty(\Omega)$ so that $0 \leq \varphi_0 \leq 1$ in $\Omega$ and
$\varphi_0 = 1$ on $\omega$. Applying (\ref{2.2}) with $w = |v_n - v|$ and
$\varphi = \varphi_0$, we have
\begin{eqnarray*}
\int_\Omega |v_n - v| \, d\mu & \leq & \int_\omega |v_n - v| \, d\mu + 2 C
|\mu|(\Omega \backslash \omega)\\
& \leq & \int_\Omega |v_n - v| \varphi_0 \, d\mu + 2 C \eps =
\int_{\Omega} |v_n - v| \varphi_0 f + \int_\Omega \nabla g \cdot \nabla\big(|v_n - v| 
\varphi_0 \big) + 2 C \eps.
\end{eqnarray*}
By (\ref{2.5}), we know that $\int_\Omega \nabla g \cdot \nabla\big(|v_n - v|
\varphi_0 \big) \to 0$ as $n \to \infty$. Since $(v_n)$ is
bounded in $L^\infty$ and $v_n \to v$ in $L^2$, we have $v_n
\rightharpoonup v$ with respect to the weak$^*$ topology of
$L^\infty$; thus, $ \int_{\Omega} |v_n - v| \varphi_0 f \to
0$. We conclude that
$\limsup_{n \to \infty}{\int_\Omega |v_n - v| \, d\mu} \leq 2 C \eps$.
Taking $\eps > 0$ arbitrarily small, (\ref{2.1}) follows.

\medskip

Given $k > 0$, we denote by $T_k : \RR \to \RR$ the truncation
operator, i.e. $T_k(s) = s$ if $s \in [-k, k]$ and $T_k(s) =
\sign{(s)} \,
k$ if $|s| > k$. Recall the following standard inequality (see e.g.\@
\cite[Lemma~1]{BrePon:03}): 

\smallskip
\begin{lemma}\label{lemma2.2}
Assume $u \in L^1\loc(\Omega)$ and $\Delta u \in \cM(\Omega)$. Then,
$T_k(u) \in H^1\loc(\Omega)$, $\forall k > 0$; moreover, given $\omega 
\subset\subset \omega' \subset\subset \Omega$, there exists $C > 0$ such that
\begin{eqnarray}\label{2.6}
\int_\omega \big|\nabla T_k(u)\big|^2 \leq k \bigg( \int_{\omega'} |\Delta u|
+ C \int_{\omega'} |u| \bigg).
\end{eqnarray}
\end{lemma}
\smallskip

Another ingredient to prove (\ref{1.3}) is our next result, which
extends Lemma~2 in \cite{BreCazMarRam:96}:

\smallskip
\begin{proposition}\label{prop2.1}
Let $\Phi : \RR \to \RR$ be a $C^1$-convex function such that $0 \leq
\Phi' \leq 1$ on $\RR$. If $u \in L^1\loc(\Omega)$ and $\Delta u \in
\cM(\Omega)$, then
\begin{equation}\label{2.7}
\Delta \Phi(u) \geq \Phi'(u) (\Delta u)\ld - (\Delta u)\lc^- \quad
\mbox{in $\cD'(\Omega)$.}
\end{equation}
\end{proposition}
\smallskip

\noindent
\textbf{Proof.} Without loss of generality, we shall assume that $\Phi
\in C^2$ and 
$\Phi''$ has compact support in $\RR$. The general case can be easily
deduced by approximation (note that since $\Phi$ is convex and $\Phi'$
 is uniformly bounded, both limits $\Phi'(\pm \infty)$ exist and are finite).
We may also assume that $u \in L^1(\Omega)$ and $\int_\Omega |\Delta
u| < \infty$.

For every $x \in \Omega$, define
$u_n(x) = \rho_n * u(x) = \int_\Omega \rho_n(x-y) u(y) \, dy$,
where $\rho_n$ is a family of radial mollifiers such that
$\supp{\rho_n} \subset B_{1/n}$. Since $\Phi'' \geq 0$
in $\RR$, we have
\begin{eqnarray*}
\Delta \Phi(u_n) = \Phi'(u_n) \Delta u_n + \Phi''(u_n) |\nabla u_n|^2
\geq \Phi'(u_n) \Delta u_n \quad \mbox{in $\Omega$.}
\end{eqnarray*}
Let $\varphi \in C_0^\infty(\Omega)$ with $\varphi \geq 0$. We
multiply both sides of the inequality above by $\varphi$ and integrate
by parts. For every $n \geq 1$ such that
$d(\supp{\varphi}, \partial\Omega) > 1/n$, we have
\begin{eqnarray*}
\int_\Omega \Phi(u_n) \Delta \varphi 
& \geq & \int_\Omega \Phi'(u_n) \varphi \, \Delta u_n\\
& = & \int_\Omega \Big\{ \rho_n * \big[ \Phi'(u_n) \varphi
  \big] \Big\} \, \Delta u \geq \int_\Omega \Big\{ \rho_n * \big[ \Phi'(u_n) \varphi
  \big] \Big\} \, (\Delta u)\ld - \int_\Omega (\rho_n * \varphi) \, (\Delta u)\lc^-.
\end{eqnarray*}
Clearly, 
\begin{equation}\label{2.8}
\int_\Omega \Phi(u_n) \Delta \varphi \to \int_\Omega \Phi(u) \Delta
\varphi \qquad \mbox{and} \qquad  \int_\Omega (\rho_n * \varphi) \,
(\Delta u)\lc^- \to  \int_\Omega \varphi \, (\Delta u)\lc^-.
\end{equation}
We now establish the following:
\smallskip

\noindent
\textbf{Claim.}
$\rho_n * \big[ \Phi'(u_n) \varphi \big] \rightharpoonup \Phi'(u)
\varphi$ in $H^1(\Omega)$.

\smallskip
In fact, since $\rho_n * \big[ \Phi'(u_n) \varphi \big] \to \Phi'(u)
\varphi$ in, say, $L^1(\Omega)$ and since $\varphi$ has compact support in
$\Omega$, it suffices to show that $\big( \Phi'(u_n) \big)$ is bounded
in $H^1\loc(\Omega)$. 
Let $M > 0$ be such that $\supp{\Phi''} \subset [-M,M]$. Then,
$$
\nabla \Phi' (u_n) = \Phi''(u_n) \nabla u_n = \Phi''(u_n) \nabla
T_M(u_n) \quad \mbox{in $\Omega$.}
$$
Let $\omega \subset\subset \omega'  \subset\subset \Omega$. For $n \geq 1$
sufficiently large, it follows from (\ref{2.6}) that
\begin{eqnarray*}
\int_\omega \big|\nabla \Phi'(u_n)\big|^2 \leq \|\Phi''\|_\infty \int_\omega \big|\nabla
T_M(u_n)\big|^2 \leq C M \bigg( \int_{\omega'} |u_n| + \int_{\omega'} |\Delta u_n| \bigg)
\leq C M \bigg( \int_\Omega |u| + \int_\Omega |\Delta u| \bigg),
\end{eqnarray*}
for some constant $C > 0$ independent of $n$.

\smallskip
In view of the previous claim, we can now apply Lemma~\ref{lemma2.1} above with
$v_n = \rho_n * \big[ \Phi'(u_n) \varphi \big]$ and $\mu = (\Delta
u)\ld$ to conclude that
\begin{equation}\label{2.9}
\int_\Omega \Big\{ \rho_n * \big[ \Phi'(u_n) \varphi
  \big] \Big\} \, (\Delta u)\ld \to \int_\Omega \Phi'(u) \varphi
\, (\Delta u)\ld.
\end{equation}
Combining (\ref{2.8}) and (\ref{2.9}) yields
\begin{eqnarray*}
\int_\Omega \Phi(u) \Delta \varphi \geq \int_\Omega \Phi'(u) \varphi
\, (\Delta u)\ld - \int_\Omega \varphi \, (\Delta u)\lc^-, \quad
\forall \varphi \in C_0^\infty(\Omega) \mbox{ with } \varphi \geq 0
\mbox{ in } \Omega,
\end{eqnarray*}
which is precisely (\ref{2.7}).

\medskip
\noindent
\textbf{Proof of (\ref{1.3}).} Let $(\Phi_n)$ be a sequence of smooth
convex functions in $\RR$ such that $\Phi_n(t) = t$ if $t \geq 0$ and
$\big| \Phi_n(t) \big| \leq 1/n$ if $t < 0$. In particular, $0 \leq
\Phi' \leq 1$ in $\RR$. It follows from the previous proposition that
\begin{eqnarray*}
\Delta \Phi_n(u) \geq \Phi_n'(u) (\Delta u)\ld - (\Delta u)\lc^- \quad
\mbox{in $\cD(\Omega)$.} 
\end{eqnarray*}
As $n \to \infty$, we get
\begin{eqnarray}\label{2.10}
\Delta u^+ \geq \chi_{[u \geq 0]} (\Delta u)\ld - (\Delta u)\lc^- \quad
\mbox{in $\cD(\Omega)$.} 
\end{eqnarray}
In particular, $\Delta u^+ \in \cM(\Omega)$. Taking the diffuse
part from both sides of (\ref{2.10}), we conclude that (\ref{1.3}) holds.

%%%%%%%%%%%%%%%%%%%%%%%%%%%%%%%%%%%%%%%%
%%%%%%%%%%%%%%%%%%%%%%%%%%%%%%%%%%%%%%%%
%%%%%%%%%%%%%%%%%%%%%%%%%%%%%%%%%%%%%%%%

\section{Proof of (\ref{1.4}) in Theorem~\ref{thm1.1}}

Identity (\ref{1.4}) relies on the following: 

\smallskip
\begin{theorem}[``Inverse'' maximum principle~\cite{DupPon:04}]\label{lemma3.1}
Let $u \in L^1\loc(\Omega)$ be such that $\Delta u \in \cM(\Omega)$. If $u
\geq 0$ a.e. in $\Omega$, then
\begin{equation}\label{3.1}
(- \Delta u)\lc \geq 0 \quad \mbox{on $\Omega$.}
\end{equation}
\end{theorem}
\smallskip

To complete the proof of Theorem~\ref{thm1.1}, we now present:

\medskip
\noindent
\textbf{Proof of (\ref{1.4}).} From the proof of (\ref{1.3}), we
already know that $\Delta u^+$ is a Radon measure on
$\Omega$. Applying the ``inverse'' maximum principle to $u^+$, we have
$(- \Delta u^+)\lc \geq 0$ on $\Omega$. Since $u^+ - u \geq 0$ a.e. in
$\Omega$, it also follows from Theorem~\ref{lemma3.1} above that $(- \Delta u^+)\lc
\geq (- \Delta u)\lc$ on $\Omega$. Thus,
\begin{eqnarray*}
(- \Delta u^+)\lc \geq (- \Delta u)\lc^+ \quad  \mbox{on $\Omega$,}
\end{eqnarray*}
which gives the ``$\geq$'' in (\ref{1.4}). The reverse inequality just
follows by taking the concentrated part from both sides of
(\ref{2.10}). In fact,
\begin{eqnarray*}
(- \Delta u^+)\lc \leq (\Delta u)\lc^- = (- \Delta u)\lc^+ \quad  \mbox{on $\Omega$.}
\end{eqnarray*}

% \appendix

% \section{}
% \label{}

% The Acknowledgements are also a un-numbered section
\section*{Acknowledgements}
The first author (H.B.) is partially sponsored
by an EC Grant through the RTN Program ``Front-Singularities'',
HPRN-CT-2002-00274. He is also a member of the Institut Universitaire
de France. 
% The second author (A.C.P.) is partially supported by CAPES,
%Brazil, under grant no. BEX1187/99-6.


\begin{thebibliography}{0}

\bibitem{Anc:79}
A. Ancona,
Une propri\'et\'e d'invariance des ensembles absorbants par
perturbation d'un op\'erateur elliptique,
Comm. Partial Differential Equations 4 (1979), 
321--337.

\bibitem{BocGalOrs:96}
L. Boccardo, T. Gallou\"et, L. Orsina,
Existence and uniqueness of entropy solutions for nonlinear
elliptic equations with measure data,
Ann. Inst. H. Poincar\'e Anal. Non Lin\'eaire 13 (1996), 
539--551.

\bibitem{BreCazMarRam:96}
H. Brezis, T. Cazenave, Y. Martel, A. Ramiandrisoa, 
Blow up for $u_t - \Delta u = g(u)$ revisited,
Adv. Differential Equations 1 (1996),
73--90.

\bibitem{BrePon:03}
H. Brezis, A.C. Ponce,
Remarks on the strong maximum principle,
Differential Integral Equations 16 (2003),
1--12.

\bibitem{BreMarPon:XX}
H. Brezis, M. Marcus, A.C. Ponce,
Nonlinear elliptic equations with measures revisited,
in preparation. 

\bibitem{DupPon:04}
L. Dupaigne, A.C. Ponce,
Singularities of positive supersolutions in elliptic PDEs,
to appear in Selecta Math. (N.S.).

\bibitem{FukSatTan:91}
M.~Fukushima, K.~Sato, S.~Taniguchi,
On the closable part of pre-Dirichlet forms and the fine supports
of underlying measures,
Osaka Math. J. 28 (1991), 
517--535.

\bibitem{Kato:72}
T. Kato,
Schr\"odinger operators with singular potentials,
Israel J. Math. 13 (1972),
135--148 (1973).



% notes:

% \bibitem{label2}
% If there is a note, it should come last:
% \bibitem{label3} \note
% \end{subbibitems}


\end{thebibliography}
\end{document}